\documentclass[12pt,a4paper,fleqn]{article}
\usepackage{a4wide,amsfonts,amsmath,latexsym,amssymb,euscript,graphicx,units,mathrsfs}

\usepackage{graphicx}
\usepackage{color}
\usepackage{amssymb}
\usepackage{amssymb}
\usepackage[T1]{fontenc}
\usepackage{latexsym}
\usepackage{xypic}
\usepackage{eufrak}
\usepackage{euscript}
\usepackage{amsfonts,amsmath}
\usepackage{verbatim}
\usepackage{fancyhdr}
\usepackage[english]{babel}
\usepackage{mathrsfs}
\usepackage{units}

\newtheorem{prop}{Proposition}[section]
\newtheorem{cor}[prop]{Corollary}
\newtheorem{corollary}[prop]{Corollary}
\newtheorem{lemme}[prop]{Lemma}
\newtheorem{lemma}[prop]{Lemma}
\newtheorem{rem}[prop]{Remark}
\newtheorem{remark}[prop]{Remark}
\newtheorem{thm}[prop]{Theorem}
\newtheorem{theorem}[prop]{Theorem}

\newtheorem{example}[prop]{Example}

\renewcommand{\geq}{\geqslant}
\def\leq{\leqslant}

\newcommand{\R}{\mathbb{R}}

\newcommand{\real}{\mathbb{R}}

\def\HH{\EuFrak H}

\def\1{{\mathbf{1}}}

\def\1{{\mathbf{1}}}
\def\0.5{{\frac{1}{2}}}

\newcommand{\fin}
{ \vspace{-0.6cm}
\begin{flushright}
\mbox{$\Box$}
\end{flushright}
\noindent }

\newcommand{\qed}{\nopagebreak\hspace*{\fill}
{\vrule width6pt height6ptdepth0pt}\par}

\begin{document}

\begin{center}
{\Large{\bf Stein's method meets Malliavin calculus:\\  a short survey with new estimates}}\\~\\
by Ivan Nourdin\footnote{Laboratoire de Probabilit{\'e}s et
Mod{\`e}les Al{\'e}atoires, Universit{\'e} Pierre et Marie Curie,
Bo{\^\i}te courrier 188, 4 Place Jussieu, 75252 Paris Cedex 5,
France, {\tt ivan.nourdin@upmc.fr}} and Giovanni
Peccati\footnote{Equipe Modal'X, Universit\'{e} Paris Ouest -- Nanterre la D\'{e}fense, 200 Avenue de la République, 92000 Nanterre, and LSTA, Universit\'{e} Paris VI, France. Email: \texttt{giovanni.peccati@gmail.com}} \\
{\it Universit\'e Paris VI and Universit\'e Paris Ouest}\\~\\
\end{center}
{\small \noindent {\bf Abstract:} We provide an overview of some recent techniques involving the Malliavin calculus of variations and the so-called ``Stein's method'' for the Gaussian approximations of probability distributions. Special attention is devoted to establishing explicit connections with the classic method of moments: in particular, we use interpolation techniques in order to deduce some new estimates for the moments of random variables belonging to a fixed Wiener chaos. As an illustration, a class of central limit theorems associated with the quadratic variation of a fractional Brownian motion is studied in detail. \\

\noindent {\bf Key words:} Central limit theorems; Fractional Brownian motion; Isonormal Gaussian processes; Malliavin calculus; Multiple integrals; Stein's method.\\

\noindent
{\bf 2000 Mathematics Subject Classification:} 60F05; 60G15; 60H05; 60H07. }
\\
\tableofcontents

\section{Introduction}
This survey deals with the powerful interaction of two probabilistic techniques, namely the {\sl Stein's method} for the normal approximation of probability distributions, and the {\sl Malliavin calculus} of variations. We will first provide an intuitive discussion of the theory, as well as an overview of the literature developed so far.

\subsection{Stein's heuristic and method} We start with an introduction to Stein's method based on moments computations. Let $N\sim\mathscr{N}(0,1)$ be a standard Gaussian random variable. It is well-known that the (integer) moments of $N$, noted $\mu_p := E(N^p)$ for $p\geq 1$, are given by: $\mu_p=0$ if $p$ is odd, and $\mu_p =(p-1)!! := p!/(2^{p/2}(p/2)!)$ if $p$ is even. A little inspection reveals that the sequence $\{\mu_p : p\geq 1\}$ is indeed completely determined by the recurrence relation:
\begin{equation}\label{introrec}
\mu_1=0, \,\,\, \mu_2 = 1,\quad \text{and}\quad \mu_p = (p-1)\times \mu_{p-2}, \quad\text{for every}\,\, p\geq 3.
\end{equation}
Now (for $p\geq 0$) introduce the notation $f_p(x) = x^p$, so that it is immediate that the relation (\ref{introrec}) can be restated as
\begin{equation}\label{introrec2}
E[N\times f_{p-1}(N)] = E[f'_{p-1}(N)], \quad \text{for every}\,\, p\geq 1.
\end{equation}
By using a standard argument based on polynomial approximations, one can easily prove that
relation (\ref{introrec2}) continues to hold if one replaces $f_p$ with a sufficiently smooth
function $f$ (e.g. any $C^1$ function with a sub-polynomial derivative will do). Now observe that
a random variable $Z$ verifying $E[Zf_{p-1}(Z)] = E[f'_{p-1}(Z)]$ for every $p\geq 1$ is
necessarily such that $E(Z^p)=\mu_p$ for every $p\geq 1$. Also, recall that
the law of a
$\mathscr{N}(0,1)$ random variable is uniquely determined by its moments.
By combining these facts with the
previous discussion, one obtains the following characterization of the (standard) normal distribution, which is universally known as ``Stein's Lemma'': {\it a random variable $Z$ has a $\mathscr{N}(0,1)$ distribution if and only if}
\begin{equation}\label{introrec3}
E[Zf(Z)-f'(Z)] =0,
\end{equation}
{\it for every smooth function $f$.} Of course, one needs to better specify the notion of ``smooth function'' -- a rigorous statement and a rigorous proof of Stein's Lemma are provided at Point 3 of Lemma \ref{stein} below.

A far-reaching idea developed by Stein (starting from the seminal paper \cite{Stein_orig}) is the following: in view of Stein's Lemma and given a generic random variable $Z$, one can measure the distance between the laws of $Z$ and $N\sim\mathscr{N}(0,1)$, by assessing the distance from zero of the quantity $E[Zf(Z) - f'(Z)]$, for every $f$ belonging to a ``sufficiently large'' class of smooth functions. Rather surprisingly, this somewhat heuristic approach to probabilistic approximations can be made rigorous by using ordinary differential equations. Indeed, one of the main findings of \cite{Stein_orig} and \cite{Stein_book} is that bounds of the following type hold in great generality:
\begin{equation}\label{introrec4}
d(Z,N) \leq C \times \sup_{f\in\mathcal{F}} |E[Zf(Z) - f'(Z)]|,
\end{equation}
where: (i) $Z$ is a generic random variable, (ii) $N\sim\mathscr{N}(0,1)$, (iii) $d(Z,N)$ indicates an appropriate distance between the laws of $Z$ and $N$ (for instance, the Kolmogorov, or the total variation distance), (iv) $\mathcal{F}$ is some appropriate class of smooth functions, and (v) $C$ is a universal constant. The case where $d$ is equal to the Kolmogorov distance, noted $d_{Kol}$, is worked out in detail in the forthcoming Section \ref{SS : steinL}: we anticipate that, in this case, one can take $C=1$, and $\mathcal{F}$ equal to the collection of all bounded Lipschitz functions with Lipschitz constant less or equal to 1.

Of course, the crucial issue in order to put Stein-type bounds into effective use, is how to assess quantities having the form of the RHS of (\ref{introrec4}). In the last thirty years, an impressive panoply of approaches has been developed in this direction: the reader is referred to the two surveys by Chen and Shao \cite{chen-shao} and Reinert \cite{Reinert_sur} for a detailed discussion of these contributions. In this paper, we shall illustrate how one can effectively estimate a quantity such as the RHS of (\ref{introrec4}), whenever the random variable $Z$ can be represented as a regular functional of a generic and possibly infinite-dimensional Gaussian field. Here, the correct notion of regularity is related to Malliavin-type operators.

\subsection{The role of Malliavin calculus}
All the definitions concerning Gaussian analysis and Malliavin calculus used in the Introduction
will be detailed in the subsequent Section \ref{S : Pre}. Let $X = \{X(h) : h\in\EuFrak{H}\}$ be
an isonormal Gaussian process over some real separable Hilbert space $\EuFrak{H}$. Suppose $Z$ is
a centered functional of $X$, such that $E(Z)=0$ and $Z$ is differentiable in the sense of
Malliavin calculus. According to the Stein-type bound (\ref{introrec4}), in order to evaluate the
distance between the law of $Z$ and the law of a Gaussian random variable $N\sim\mathscr{N}(0,1)$,
one must be able to assess the distance between the two quantities $E[Zf(Z)]$ and $E[f'(Z)]$. The
main idea developed in \cite{NP-PTRF}, and later in the references \cite{exact, Noupecrei2, Noupecrei3, multivariate}, is that the needed estimate can be realized by using the following consequence of the {\sl integration by parts formula} of Malliavin calculus: for every $f$ sufficiently smooth (see Section \ref{SS :Mall} for a more precise statement),
\begin{equation}\label{bestformulaever}
E[Zf(Z)] = E[f'(Z)\langle DZ, -DL^{-1}Z \rangle_\EuFrak{H}],
\end{equation}
where $D$ is the Malliavin derivative operator, $L^{-1}$ is the pseudo-inverse of the Ornstein-Uhlenbeck generator, and $\langle\cdot , \cdot \rangle_\EuFrak{H}$ is the inner product of $\EuFrak{H}$. It follows from (\ref{bestformulaever}) that, if the derivative $f'$ is bounded, then the distance between $E[Zf(Z)]$ and $E[f'(Z)]$ is controlled by the $L^1(\Omega)$-norm of the random variable $1 - \langle DZ, -DL^{-1}Z \rangle_\EuFrak{H}$. For instance, in the case of the Kolmogorov distance, one obtains that, for every centered and Malliavin differentiable random variable $Z$,
\begin{equation}\label{nicebound}
d_{Kol}(Z,N) \leq E|1-\langle DZ, -DL^{-1}Z \rangle_\EuFrak{H}|.
\end{equation}
We will see in Section \ref{SS : fourth} that, in the particular case where $Z=I_q(f)$ is a multiple Wiener-It\^{o} integral of order $q\geq 2$ (that is, $Z$ is an element of the $q$th Wiener chaos of $X$) with unit variance, relation (\ref{nicebound}) yields the neat estimate
\begin{equation}\label{nicebound2}
d_{Kol}(Z,N) \leq \sqrt{\frac{q-1}{3q}\times |E(Z^4)-3|}.
\end{equation}
Note that $E(Z^4)-3$ is just the fourth cumulant of $Z$, and that the fourth cumulant of $N$ equals zero. We will also show that the combination of (\ref{nicebound}) and (\ref{nicebound2}) allows to recover (and refine) several characterizations of CLTs on a fixed Wiener chaos -- as recently proved in \cite{NO} and \cite{nunugio}.

\subsection{Beyond the method of moments} The estimate (\ref{nicebound2}), specially when combined with the findings of \cite{multivariate} and \cite{PTu04} (see Section \ref{S : MULTI}), can be seen as a drastic simplification of the so-called ``method of moments and cumulants'' (see Major \cite{Major} for a classic discussion of this method in the framework of Gaussian analysis). Indeed, such a relation implies that, if $\{Z_n : n\geq 1\}$ is a sequence of random variables with unit variance belonging to a fixed Wiener chaos, then, in order to prove that $Z_n$ converges in law to $N\sim\mathscr{N}(0,1)$, it is sufficient to show that $E(Z_n^4)$ converges to $E(N^4)=3$. Again by virtue of (\ref{nicebound2}), one also has that the rate of convergence of $E(Z_n^4)$ to 3 determines the ``global'' rate convergence in the Kolmogorov distance. In order to further characterize the connections between our techniques and moments computations, in Proposition \ref{P : moments} we will deduce some new estimates, implying that (for $Z$ with unit variance and belonging to a fixed Wiener chaos), for every integer $k\geq 3$ the quantity $|E[Z^k] - E[N^k]|$ is controlled (up to an explicit universal multiplicative constant) by the square root of $|E[Z^4] - E[N^4]|$. This result is obtained by means of an interpolation technique, recently used in \cite{Noupecrei3} and originally introduced by Talagrand -- see e.g. \cite{tala}.

\subsection{An overview of the existing literature} The present survey is mostly based on the three references \cite{NP-PTRF}, \cite{Noupecrei3} and \cite{multivariate}, dealing with upper bounds in the one-dimensional and multi-dimensional approximations of regular functionals of general Gaussian fields (strictly speaking, the papers \cite{NP-PTRF} and \cite{Noupecrei3} also contain results on non-normal approximations, related e.g. to the Gamma law). However, since the appearance of \cite{NP-PTRF}, several related works have been written, which we shall now shortly describe.
\begin{itemize}
\item[-] Our paper \cite{exact} is again based on Stein's method and Malliavin calculus, and deals with the problem of determining optimal rates of convergence. Some results bear connections with one-term Edgeworth expansions.
\item[-] The paper \cite{breton-nourdin}, by Breton and Nourdin, completes the analysis initiated in \cite[Section 4]{NP-PTRF}, concerning the obtention of Berry-Ess\'{e}en bounds associated with the so-called Breuer-Major limit theorems (see \cite{BrMa}). The case of non-Gaussian limit laws (of the Rosenblatt type) is also analyzed.
\item[-] In \cite{Noupecrei2}, by Nourdin, Peccati and Reinert, one can find an application of Stein's method and Malliavin calculus to the derivation of second order Poincar\'{e} inequalities on Wiener space. This also refines the CLTs on Wiener chaos proved in \cite{NO} and \cite{nunugio}.
\item[-] One should also mention our paper \cite{noncentral}, where we deduce a characterization of non-central limit theorems (associated with Gamma laws) on Wiener chaos. The main findings of \cite{noncentral} are refined in \cite{NP-PTRF} and \cite{Noupecrei3}, again by means of Stein's method.
\item[-] The work \cite{NV}, by Nourdin and Viens, contains an application of (\ref{bestformulaever}) to the estimate of densities and tail probabilities associated with functionals of Gaussian processes, like for instance quadratic functionals or suprema of continuous-time Gaussian processes on a finite interval.
\item[-] The findings of \cite{NV} have been further refined by Viens in \cite{Viens}, where one can also find some applications to polymer fluctuation exponents.
\item[-] The paper \cite{breton-nourdin-peccati}, by Breton, Nourdin and Peccati, contains some statistical applications of the results of \cite{NV}, to the construction of confidence intervals for the Hurst parameter of a fractional Brownian motion.
\item[-] Reference \cite{BeTaNou}, by Bercu, Nourdin and Taqqu, contains some applications of the results of \cite{NP-PTRF} to almost sure CLTs.
\item[-] In \cite{Noupecrei1}, by Nourdin, Peccati and Reinert, one can find an extension of the ideas introduced in \cite{NP-PTRF} to the framework of functionals of Rademacher sequences. To this end, one must use a discrete version of Malliavin calculus (see Privault \cite{Privault}).
\item[-] Reference \cite{PecSoleTaqUtz}, by Peccati, Sol\'{e}, Taqqu and Utzet, concerns a combination of Stein's method with a version of Malliavin calculus on the Poisson space (as developed by Nualart and Vives in \cite{NuVives}).
\item[-] Reference \cite{Noupecrei3}, by Nourdin, Peccati and Reinert, contains an application of Stein's method, Malliavin calculus and the ``Lindeberg invariance principle'', to the study of universality results for sequences of homogenous sums associated with general collections of independent random variables.
\end{itemize}

\section{Preliminaries}\label{S : Pre}
We shall now present the basic elements of Gaussian analysis and Malliavin calculus that are used in this paper. The reader
is referred to the monograph by  Nualart \cite{nualartbook} for any unexplained definition or result.
\subsection{Isonormal Gaussian processes}
Let $\EuFrak H$ be a real separable Hilbert space. For any $q\geq 1$, we denote by $\EuFrak H^{\otimes q}$ the $q$th tensor product of $\EuFrak H$, and by
$\EuFrak H^{\odot q}$ the associated $q$th {\sl symmetric} tensor product; plainly, $\EuFrak H^{\otimes 1}=\EuFrak H^{\odot 1}= \EuFrak H$.

We write $X=\{X(h),h\in \EuFrak H\}$ to indicate
an \textsl{isonormal Gaussian process} over
$\EuFrak H$. This means that $X$ is a centered Gaussian family, defined on some probability space $(\Omega ,\mathcal{F},P)$, and such that $E\left[ X(g)X(h)\right] =\langle g,h\rangle _{\EuFrak H}$ for every $g,h\in \EuFrak H$. Without loss of generality, we also assume that
$\mathcal{F}$ is generated by $X$.

The concept of an isonormal Gaussian process dates back to Dudley's paper \cite{Dudley}. As shown in the forthcoming five examples, this general notion may be used to encode the structure of many remarkable Gaussian families.

\begin{example}[Euclidean spaces] {\rm Fix an integer $d\geq 1$, set $%
\mathfrak{H}=\mathbb{R}^{d}$ and let $\left( e_{1},...,e_{d}\right) $ be an
orthonormal basis of $\mathbb{R}^{d}$ (with respect to the usual Euclidean
inner product). Let $\left( Z_{1},...,Z_{d}\right) $ be a Gaussian vector
whose components are i.i.d. $N\left( 0,1\right) $. For every $%
h=\sum_{j=1}^{d}c_{j}e_{j}$ (where the $c_{j}$ are real and uniquely
defined), set $X\left( h\right) =\sum_{j=1}^{d}c_{j}Z_{j}$ and define $%
X=\left\{ X\left( h\right) :h\in \mathbb{R}^{d}\right\} $. Then, $X$ is an
isonormal Gaussian process over $\mathbb{R}^{d}$ endowed with its canonical
inner product.}
\end{example}

\begin{example}[Gaussian measures]{\rm Let $\left( A,\mathcal{A},\nu
\right) $ be a measure space, where $\nu $ is positive, $\sigma$-finite and
non-atomic. Recall that a (real) \textsl{Gaussian random measure} over $( A,\mathcal{A})$, with control $\nu$, is a centered Gaussian family of the type
$$G=\left\{ G\left(
B\right) :B\in \mathcal{A},\, \nu(B)<\infty \right\},$$
satisfying the relation: for every $B,C\in \mathcal{A}$ of finite $\nu$-measure, $E[G(B)G(C)] = \nu(B \cap C)$.
Now consider the Hilbert space $\mathfrak{H}=L^{2}\left( A,\mathcal{A},\nu \right) $, with inner product $\langle h,h'\rangle_\mathfrak{H} \!=\!\! \int_A\! h(a)h'(a)\nu(da)$. For every $h\in
\mathfrak{H}$, define $X\left( h\right) =\int_A h(a)G(da) $ to be
the Wiener-It\^{o} integral of $h$ with respect to $G$. Then, $X=\left\{ X\left(
h\right) :h\in L^{2}\left( Z,\mathcal{Z},\nu \right) \right\} $ defines a centered Gaussian
family with covariance given by $E[X(h)X(h')] = \langle h,h'\rangle_\mathfrak{H}$, thus yielding that $X$ is an
isonormal Gaussian process over $L^{2}\left( A,\mathcal{A},\nu \right)$. For instance, by setting $A=[0,+\infty)$ and $\nu$ equal to the Lebesgue measure, one obtains that the process $W_t = G([0,t))$, $t\geq 0$, is a standard Brownian motion started from zero (of course, in order to meet the usual definition of a Brownian motion, one has also to select a continuous version of $W$), and $X$ coincides with the $L^2(\Omega)$-closed linear Gaussian space generated by $W$.}
\end{example}

\begin{example} [Isonormal spaces derived from covariances]\label{Exemple : cov}{\rm Let $%
Y=\left\{ Y_{t}:t\geq 0\right\} $ be a real-valued centered Gaussian\
process indexed by the positive axis, and set $R\left( s,t\right) ={E}%
\left[ Y_{s}Y_{t}\right] $ to be the covariance function of $Y$. One
can embed $Y$ into some isonormal Gaussian process as follows: (i) define $%
\mathscr{E}$\ as the collection of all finite linear combinations of
indicator functions of the type $\mathbf{1}_{\left[ 0,t\right] }$, $t\geq 0$%
; (ii) define $\mathfrak{H}=\mathfrak{H}_{R}$ to be the Hilbert space given
by the closure of $\mathscr{E}$ with respect to the inner product%
\begin{equation*}
\langle f,h\rangle _{R}:=\sum_{i,j}a_{i}c_{j}R\left( s_{i},t_{j}\right) \text{%
,}
\end{equation*}%
where $f=\sum_{i}a_{i}\mathbf{1}_{\left[ 0,s_{i}\right] }$ and $%
h=\sum_{j}c_{j}\mathbf{1}_{\left[ 0,t_{j}\right] }\ $are two generic
elements of $\mathscr{E}$; (iii) for $h=\sum_{j}c_{j}\mathbf{1}_{\left[
0,t_{j}\right] }\in \mathcal{E}$, set $X\left( h\right)
=\sum_{j}c_{j}Y_{t_{j}}$; (iv) for $h\in \mathfrak{H}_{R}$, set $X\left(
h\right) $\ to be the $L^{2}\left( {P}\right) $ limit\ of any
sequence of the type $X\left( h_{n}\right) $, where $\left\{ h_{n}\right\}
\subset \mathscr{E}$ converges to $h$ in $\mathfrak{H}_{R}$. Note that such
a sequence $\left\{ h_{n}\right\} $ necessarily exists and may not be unique
(however, the definition of $X\left( h\right) $ does not depend on the
choice of the sequence $\left\{ h_{n}\right\} $). Then, by construction, the
Gaussian space $\left\{ X\left( h\right) :h\in \mathfrak{H}_R\right\} $ is an
isonormal Gaussian process over $\mathfrak{H}_{R}$. See Janson \cite[Ch. 1]%
{Janson} or Nualart \cite{nualartbook}, as well as the forthcoming Section \ref{BreuMaj}, for more details on this construction.
}\end{example}

\begin{example}[Even functions and symmetric measures]{\rm
Other classic examples of isonormal Gaussian processes (see e.g., \cite%
{ChaSlud, GiSu, Major}) are given by objects of
the type $$X_{\beta }=\left\{ X_{\beta }\left( \psi \right) :\psi \in
\mathfrak{H}_{\sl{E},\beta }\right\},$$ where $\beta $ is a real
non-atomic symmetric measure on $\left( -\pi ,\pi \right] $ (that is, $\beta
\left( dx\right) =\beta \left( -dx\right) $), and
\begin{equation}
\mathfrak{H}_{\sl{E},\beta }=L_{\sl{E}}^{2}\left( \left( -\pi ,\pi %
\right] ,d\beta \right)  \label{evenN1}
\end{equation}%
stands for the collection of all \textsl{real} linear combinations of
complex-valued \textsl{even} functions that are square-integrable with
respect to $\beta $ (recall that a function $\psi $ is even if $\overline{%
\psi \left( x\right) }=\psi \left( -x\right) $). The class $\mathfrak{H}_{%
\sl{E},\beta }$ is a real Hilbert space, endowed with the inner
product
\begin{equation}
\langle \psi _{1},\psi _{2}\rangle _{\beta }=\int_{-\pi }^{\pi }\psi
_{1}\left( x\right) \psi _{2}\left( -x\right) \beta \left( dx\right) \in
\mathbb{R}.  \label{evenN2}
\end{equation}
This type of construction is used in the spectral theory of time series.}
\end{example}

\begin{example}[Gaussian Free Fields] {\rm Let $d\geq 2$ and let $D$ be a domain in $\mathbb{R}^d$. Denote by $H_s(D)$ the space of real-valued continuous and continuously differentiable functions on $\mathbb{R}^d$ that are supported on a compact subset of $D$ (note that this implies that the first derivatives of the elements of $H_s(D)$ are square-integrable with respect to the Lebesgue measure). Write $H(D)$ in order to indicate real Hilbert space obtained as the closure of $H_s(D)$ with respect to the inner product $\langle f,g \rangle = \int_{\mathbb{R}^d} \nabla f(x)\cdot \nabla g(x) dx$, where $\nabla$ is the gradient. An isonormal Gaussian process of the type $X = \{X(h) : h\in H(D)\}$ is called a \textsl{Gaussian Free Field} (GFF). The reader is referred to the survey by Sheffield \cite{SheffieldGFF} for a discussion of the emergence of GFFs in several areas of modern probability. See e.g. Rider and Vir\'{a}g \cite{RiderVirag} for a connection with the ``circular law'' for Gaussian non-Hermitian random matrices.}
\end{example}

\begin{rem}{\rm
An isonormal Gaussian process is simply an isomorphism between a centered $L^2(\Omega)$-closed
linear Gaussian space and a real separable Hilbert space $\mathfrak{H}$.
Now, fix a generic centered $L^2(\Omega)$-closed linear Gaussian space, say $\mathcal{G}$. Since $\mathcal{G}$ is itself a real separable Hilbert space (with respect to the usual $L^2(\Omega)$ inner product) it follows that $\mathcal{G}$ can always be (trivially) represented as an isonormal Gaussian process, by setting $\mathfrak{H}=\mathcal{G}$. Plainly, the subtlety in the use of isonormal Gaussian processes is that one has to select an isomorphism that is well-adapted to the specific problem one wants to tackle.
}
\end{rem}

\subsection{Chaos, hypercontractivity and products}
We now fix a generic isonormal Gaussian process $X=\{X(h),h\in \EuFrak H\}$, defined on some space $(\Omega,\mathcal{F},P)$ such that $\sigma(X)=\mathcal{F}$.

\smallskip

\noindent {\bf Wiener chaos.} For every $q\geq 1$, we write $\mathcal{H}_{q}$ in order to indicate the $q$th {\sl Wiener chaos} of $X$.
We recall that $\mathcal{H}_{q}$ is the closed linear subspace of $L^2(\Omega ,\mathcal{F},P)$
generated by the random variables of the type $H_{q}(X(h))$, where $h\in \EuFrak H$ is such that $\left\|
h\right\| _{\EuFrak H}=1$, and $H_{q}$ stands for the $q$th {\sl Hermite polynomial},
defined as
\begin{equation}\label{hermite}
H_q(x) = (-1)^q e^{\frac{x^2}{2}}
 \frac{d^q}{dx^q}  e^{-\frac{x^2}{2}},\quad x\in\R, \quad q\geq 1.
 \end{equation}
We also use the convention $\mathcal{H}_{0} = \mathbb{R}$. For
any $q\geq 1$, the mapping
\begin{equation}\label{hermite-link}
I_{q}(h^{\otimes q})=q!H_{q}(X(h))
\end{equation}
can be extended to a
linear isometry between the symmetric tensor product $\EuFrak H^{\odot q}$
equipped with the modified norm $\sqrt{q!}\left\| \cdot \right\| _{\EuFrak H^{\otimes q}}$ and the $q$th Wiener chaos $\mathcal{H}_{q}$. For $q=0$, we write $I_{0}(c)=c$, $c\in\mathbb{R}$.

\begin{rem}\label{Rem : Sym}{\rm When $\mathfrak{H}=L^{2}\left( A,\mathcal{A},\nu \right)$, the symmetric tensor product $\EuFrak H^{\odot q}$ can be identified with the Hilbert space $L_s^{2}\left( A^q,\mathcal{A}^q,\nu^q \right)$, which is defined as the collection of all symmetric functions on $A^q$ that are square-integrable with respect to $\nu^q$. In this case, it is well-known that the random variable $I_q(h)$, $h\in \EuFrak H^{\odot q}$, coincides with the (multiple) {\sl Wiener-It\^{o} integral}, of order $q$, of $h$ with respect to the Gaussian measure $B \mapsto X({\bf 1}_B)$, where $B\in\mathcal{A}$ has finite $\nu$-measure. See \cite[Chapter 1]{nualartbook} for more details on this point.
}
\end{rem}

\noindent {\bf Hypercontractivity.} Random variables living in a fixed Wiener chaos are hypercontractive. More precisely,
assume that $Z$ belongs
to the $q$th Wiener chaos $\mathcal{H}_{q}$ ($q\geq 1$). Then, $Z$ has a finite variance by construction and, for all $p\in [ 2,+\infty)$, one has the following estimate (see \cite[Th. 5.10]{Janson} for a proof):
\begin{equation}\label{hypercontractivity}
E\big(|Z|^p\big)\leq (p-1)^{pq/2}E\big(Z^2\big)^{p/2}.
\end{equation}
In particular, if $E(Z^2) =1$, one has that $E\big(|Z|^p\big)\leq (p-1)^{pq/2}$. For future use,
we also observe that, for every $q\geq 1$, the mapping $p \mapsto (p-1)^{pq/2}$ is strictly
increasing on $[2,+\infty)$.

\medskip

\noindent {\bf Chaotic decompositions.} It is well-known ({\sl Wiener chaos decomposition}) that
the space $L^2(\Omega ,\mathcal{F},P)$
can be decomposed into the infinite orthogonal sum of the spaces $\mathcal{H}_{q}$. It follows
that any square-integrable random variable
$Z\in L^2(\Omega ,\mathcal{F},P)$ admits the following chaotic expansion
\begin{equation}
Z=\sum_{q=0}^{\infty }I_{q}(f_{q}),  \label{E}
\end{equation}
where $f_{0}=E(Z)$, and the kernels $f_{q}\in \EuFrak H^{\odot q}$, $q\geq 1$, are
uniquely determined. For every $q\geq 0$, we also denote by $J_{q}$ the
orthogonal projection operator on $\mathcal{H}_q$. In particular, if
$Z\in L^2(\Omega ,\mathcal{F},P)$ is as in (\ref{E}), then
$J_{q}(Z)=I_{q}(f_{q})$ for every $q\geq 0$.

\medskip

\noindent {\bf Contractions.} Let $\{e_{k},\,k\geq 1\}$ be a complete orthonormal system in $\EuFrak H$.
Given $f\in \EuFrak H^{\odot p}$ and $g\in \EuFrak H^{\odot q}$, for every
$r=0,\ldots ,p\wedge q$, the {\sl contraction} of $f$ and $g$ of order $r$
is the element of $\EuFrak H^{\otimes (p+q-2r)}$ defined by
\begin{equation}
f\otimes _{r}g=\sum_{i_{1},\ldots ,i_{r}=1}^{\infty }\langle
f,e_{i_{1}}\otimes \ldots \otimes e_{i_{r}}\rangle _{\EuFrak H^{\otimes
r}}\otimes \langle g,e_{i_{1}}\otimes \ldots \otimes e_{i_{r}}
\rangle_{\EuFrak H^{\otimes r}}.  \label{v2}
\end{equation}
Notice that $f\otimes _{r}g$ is not necessarily symmetric: we denote its
symmetrization by $f\widetilde{\otimes }_{r}g\in \EuFrak H^{\odot (p+q-2r)}$.
Moreover, $f\otimes _{0}g=f\otimes g$ equals the tensor product of $f$ and
$g$ while, for $p=q$, one has that $f\otimes _{q}g=\langle f,g\rangle _{\EuFrak H^{\otimes q}}$.
In the particular case where $\EuFrak H=L^2(A,\mathcal{A},\nu )$, one has that $\EuFrak H^{\odot q}=L_{s}^{2}(A^{q},
\mathcal{A}^q,\nu ^q)$ (see Remark \ref{Rem : Sym}) and the contraction in (\ref{v2}) can be written in integral form as
\begin{eqnarray*}
(f\otimes _{r}g)(t_1,\ldots,t_{p+q-2r})
&=&\int_{A^{r}}f(t_{1},\ldots ,t_{p-r},s_{1},\ldots ,s_{r}) \\
&&\times \,g(t_{p-r+1},\ldots ,t_{p+q-2r},s_{1},\ldots ,s_{r})d\nu
(s_{1})\ldots d\nu (s_{r}).
\end{eqnarray*}

\noindent {\bf Multiplication.}
The following multiplication formula is well-known: if $f\in \EuFrak
H^{\odot p}$ and $g\in \EuFrak
H^{\odot q}$, then
\begin{eqnarray}\label{multiplication}
I_p(f) I_q(g) = \sum_{r=0}^{p \wedge q} r! {p \choose r}{ q \choose r} I_{p+q-2r} (f\widetilde{\otimes}_{r}g).
\end{eqnarray}
Note that (\ref{multiplication})
gives an immediate proof of the fact that multiple Wiener-It\^{o} integrals have finite moments of every order.

\subsection{The language of Malliavin calculus}\label{SS :Mall}

We now introduce some basic elements of the Malliavin calculus with respect
to the isonormal Gaussian process $X$.

\medskip

\noindent {\bf Malliavin derivatives.} Let $\mathcal{S}$ be the set of all
{\sl cylindrical random variables} of
the type
\begin{equation}
Z=g\left( X(\phi _{1}),\ldots ,X(\phi _{n})\right) ,  \label{v3}
\end{equation}
where $n\geq 1$, $g:\mathbb{R}^{n}\rightarrow \mathbb{R}$ is an infinitely
differentiable function with compact support and $\phi _{i}\in \EuFrak H$.
The Malliavin derivative of $Z$ with respect to $X$ is the element of
$L^2(\Omega ,\EuFrak H)$ defined as
\begin{equation*}
DZ\;=\;\sum_{i=1}^{n}\frac{\partial g}{\partial x_{i}}\left( X(\phi
_{1}),\ldots ,X(\phi _{n})\right) \phi _{i}.
\end{equation*}
By iteration, one can
define the $m$th derivative $D^{m}Z$, which is an element of $L^2(\Omega ,\EuFrak H^{\odot m})$,
for every $m\geq 2$.
For $m\geq 1$ and $p\geq 1$, ${\mathbb{D}}^{m,p}$ denotes the closure of
$\mathcal{S}$ with respect to the norm $\Vert \cdot \Vert _{m,p}$, defined by
the relation
\begin{equation*}
\Vert Z\Vert _{m,p}^{p}\;=\;E\left[ |Z|^{p}\right] +\sum_{i=1}^{m}E\left(
\Vert D^{i}Z\Vert _{\EuFrak H^{\otimes i}}^{p}\right) .
\end{equation*}

\smallskip

\noindent {\bf The chain rule.} The Malliavin derivative $D$ verifies the following chain rule. If
$\varphi :\mathbb{R}^{d}\rightarrow \mathbb{R}$ is continuously
differentiable with bounded partial derivatives and if $Z=(Z_{1},\ldots
,Z_{d})$ is a vector of elements of ${\mathbb{D}}^{1,2}$, then $\varphi
(Z)\in {\mathbb{D}}^{1,2}$ and
\begin{equation}\label{chainrule}
D\,\varphi (Z)=\sum_{i=1}^{d}\frac{\partial \varphi }{\partial x_{i}}(Z)DZ_{i}.
\end{equation}
A careful application e.g. of the multiplication formula (\ref{multiplication}) shows that (\ref{chainrule}) continues to hold whenever the function $\varphi$ is a polynomial in $d$ variables. Note also that a random variable $Z$ as in (\ref{E}) is in ${\mathbb{D}}^{1,2}$ if and only if
$\sum_{q=1}^{\infty }q\|J_q(Z)\|^2_{L^2(\Omega)}<\infty$
and, in this case, $E\left( \Vert DZ\Vert _{\EuFrak H}^{2}\right)
=\sum_{q=1}^{\infty }q\|J_q(Z)\|^2_{L^2(\Omega)}$. If $\EuFrak H=
{L}^{2}(A,\mathcal{A},\nu )$ (with $\nu $ non-atomic), then the
derivative of a random variable $Z$ as in (\ref{E}) can be identified with
the element of $L^2(A\times \Omega )$ given by
\begin{equation}
D_{x}Z=\sum_{q=1}^{\infty }qI_{q-1}\left( f_{q}(\cdot ,x)\right) ,\quad x\in
A.  \label{dtf}
\end{equation}

\noindent {\bf The divergence operator.} We denote by $\delta $ the adjoint of the operator $D$, also called the
\textsl{divergence operator}. A random element $u\in L^2(\Omega ,\EuFrak H)$
belongs to the domain of $\delta $, noted $\mathrm{Dom}\delta $, if and
only if it verifies
$|E\langle DZ,u\rangle _{\EuFrak H}|\leq c_{u}\,\Vert Z\Vert _{L^2(\Omega)}$
for any $Z\in \mathbb{D}^{1,2}$, where $c_{u}$ is a constant depending only
on $u$. If $u\in \mathrm{Dom}\delta $, then the random variable $\delta (u)$
is defined by the duality relationship
\begin{equation}
E(Z\delta (u))=E\big(\langle DZ,u\rangle _{\EuFrak H}\big),  \label{ipp}
\end{equation}
which holds for every $Z\in {\mathbb{D}}^{1,2}$.

\smallskip

\noindent {\bf Ornstein-Uhlenbeck operators.} The operator $L$, known as the {\sl generator of the Ornstein-Uhlenbeck semigroup}, is defined as
$L=\sum_{q=0}^{\infty }-qJ_{q}$. The domain of $L$ is
\begin{equation*}
\mathrm{Dom}L=\{Z\in L^2(\Omega ):\sum_{q=1}^{\infty }q^{2}\left\|
J_{q}(Z)\right\| _{L^2(\Omega )}^{2}<\infty \}=\mathbb{D}^{2,2}\text{.}
\end{equation*}
There is an important relation between the operators $D$, $\delta $ and $L$
(see e.g. \cite[Proposition 1.4.3]{nualartbook}): a random variable $Z$ belongs to
$\mathbb{D}^{2,2}$ if and only if $Z\in \mathrm{Dom}\left( \delta D\right) $
(i.e. $Z\in {\mathbb{D}}^{1,2}$ and $DZ\in \mathrm{Dom}\delta $) and, in
this case,
\begin{equation}
\delta DZ=-LZ.  \label{k1}
\end{equation}
For any $Z \in L^2(\Omega )$, we define $L^{-1}Z =\sum_{q=1}^{\infty }-\frac{1}{q} J_{q}(Z)$. The operator $L^{-1}$ is called the
{\sl pseudo-inverse} of $L$. For any $Z \in L^2(\Omega )$, we have that $L^{-1} Z\in  \mathrm{Dom}L$,
and
\begin{equation}\label{Lmoins1}
LL^{-1} Z = Z - E(Z).
\end{equation}

\noindent {\bf An important string of identities.} Finally, let us mention a chain of identities playing a crucial role in the sequel.
Let $f:\R\to\R$ be a $C^1$ function with bounded derivative, and let
$F,Z\in \mathbb{D}^{1,2}$. Assume moreover that $E(Z)=0$.
By using successively (\ref{Lmoins1}), (\ref{k1}) and (\ref{chainrule}), one deduces that
\begin{eqnarray}
E\big(Zf(F)\big)&=&E\big(LL^{-1}Z\times f(F)\big)=E\big(\delta D(-L^{-1}Z)\times f(F)\big) \notag \\
&=& E\big(\langle Df(F),-DL^{-1}Z\rangle_\HH\big)\notag\\
&=&E\big(f'(F)\langle DF,-DL^{-1}Z\rangle_\HH\big).\label{ivangioformula}
\end{eqnarray}
We will shortly see that the fact $E\big(Zf(F)\big)=E\big(f'(F)\langle DF,-DL^{-1}Z\rangle_\HH\big)$ constitutes a fundamental element in the connection between Malliavin calculus and Stein's method.

\section{One-dimensional approximations}
\subsection{Stein's lemma for normal approximations}\label{SS : steinL}
Originally introduced in the path-breaking paper \cite{Stein_orig}, and then further developed in the monograph \cite{Stein_book}, {\sl Stein's method} can be roughly described as a collection of probabilistic techniques, allowing to characterize the approximation of probability distributions by means of differential operators. As already pointed out in the Introduction, the two surveys \cite{chen-shao} and \cite{Reinert_sur} provide a valuable introduction to this very active area of modern probability.
In this section, we are mainly interested in the use of Stein's method for the {\sl normal approximation} of the laws of real-valued random variables, where the approximation is performed with respect to the {\sl Kolmogorov distance}. We recall that
the Kolmogorov distance between the laws of two real-valued random variables $Y$ and $Z$ is defined by
\[
d_{Kol}(Y,Z)=\sup_{z\in\R} \big| P(Y\leq z)-P(Z \leq z)\big|.
\]

The reader is referred to \cite{NP-PTRF} for several extensions of the results discussed in this survey to other distances between probability measures, such as e.g. the {\sl total variation distance}, or the {\sl Wasserstein distance}. The following statement, containing all the elements of Stein's method that are needed for our discussion, can be traced back to Stein's original contribution \cite{Stein_orig}.

\begin{lemme}\label{stein}
Let $N\sim \mathscr{N}(0,1)$ be a standard Gaussian random variable.
\begin{itemize}
\item[\rm 1.] Fix $z\in\mathbb{R}$, and define
$f_z:\mathbb{R}\to\mathbb{R}$ as
\begin{equation}\label{fz}
f_z(x)=e^{\frac{x^2}{2}}\int_{-\infty}^x \big({\bf 1}_{(-\infty,z]}(a)-P(N\leq z)\big)
e^{-\frac{a^2}{2}}da,\quad x\in \R.
\end{equation}
Then, $f_z$ is continuous on $\mathbb{R}$, bounded by $\sqrt{2\pi}/4$, differentiable on
$\mathbb{R}\setminus\{z\}$, and verifies moreover
\begin{equation}\label{stein-equation}
f_z'(x)-xf_z(x)={\bf 1}_{(-\infty,z]}(x)-P(N\leq z)
\quad \mbox{for all }x\in\mathbb{R}\setminus\{z\}.
\end{equation}
One has also that $f_z$ is Lipschitz, with Lipschitz constant less or equal to 1.
\item[\rm 2.] Let $Z$ be a generic random variable. Then,
\begin{equation}\label{StCrudeBound}
d_{Kol}(Z,N) \leq \sup_f |E[Zf(Z) - f'(Z)]|,
\end{equation}
where the supremum is taken over the class of all Lipschitz functions that are bounded by $\sqrt{2\pi}/4$ and whose Lipschitz constant is less or equal to 1.
\item[\rm 3.] Let $Z$ be a generic random variable. Then, $Z \sim \mathscr{N}(0,1)$ if and only if $E[Zf(Z) - f'(Z)]=0$ for every continuous and piecewise differentiable function $f$ verifying the relation $E|f'(N)|<\infty$.
\end{itemize}
\end{lemme}
{\it Proof}:
({\it Point 1}) We shall only prove that $f_z$ is Lipschitz and we will evaluate its constant (the proof of the remaining properties is left to the reader).  We have, for $x\geq 0$, $x\neq z$:
\begin{eqnarray*}
\big|f'_z(x)\big|&=&\left|{\bf 1}_{(-\infty,z]}(x)-P(N\leq z)
+ xe^{\frac{x^2}{2}}\int_{-\infty}^x \big({\bf 1}_{(-\infty,z]}(a)-P(N\leq z)\big)
e^{-\frac{a^2}{2}}da\right|\\
&\underset{(*)}{=}&\left|{\bf 1}_{(-\infty,z]}(x)-P(N\leq z)
- xe^{\frac{x^2}{2}}\int_x^{+\infty}\big({\bf 1}_{(-\infty,z]}(a)-P(N\leq z)\big)
e^{-\frac{a^2}{2}}da\right|\\
&\leq&\big\| {\bf 1}_{(-\infty,z]}(\cdot)-P(N\leq z)\big\|_\infty
\left(1+ xe^{\frac{x^2}{2}}\int_x^{+\infty}
e^{-\frac{a^2}{2}}da\right)\\
&\leq& 1+ e^{\frac{x^2}{2}}\int_x^{+\infty}
ae^{-\frac{a^2}{2}}da=2.
\end{eqnarray*}
Observe that identity $(*)$ holds since
\[
0=E\big({\bf 1}_{(-\infty,z]}(N)-P(N\leq z)\big)
=\frac{1}{\sqrt{2\pi}}\int_{-\infty}^{+\infty}
\big({\bf 1}_{(-\infty,z]}(a)-P(N\leq z)\big)
e^{-\frac{a^2}{2}}da.
\]
For $x\leq 0$, $x\neq z$, we can write
\begin{eqnarray*}
\big|f'_z(x)\big|&=&\left|{\bf 1}_{(-\infty,z]}(x)-P(N\leq z)
+ xe^{\frac{x^2}{2}}\int_{-\infty}^x \big({\bf 1}_{(-\infty,z]}(a)-P(N\leq z)\big)
e^{-\frac{a^2}{2}}da\right|\\
&\leq&\big\| {\bf 1}_{(-\infty,z]}(\cdot)-P(N\leq z)\big\|_\infty
\left(1+ |x|e^{\frac{x^2}{2}}\int_{-\infty}^{x}
e^{-\frac{a^2}{2}}da\right)\\
&\leq& 1+ e^{\frac{x^2}{2}}\int_{-\infty}^x
|a|e^{-\frac{a^2}{2}}da=2.
\end{eqnarray*}
Hence, we have shown that $f_z$ is Lipschitz with Lipschitz constant bounded by 2.
For the announced refinement (that is, the constant is bounded by 1),
we refer the reader to Chen and Shao \cite[Lemma 2.2]{chen-shao}.

\smallskip

\noindent ({\it Point 2}) Take expectations on both sides of (\ref{stein-equation}) with respect to the law of $Z$. Then, take the supremum over all $z\in\R$, and exploit the properties of $f_z$ proved at Point 1.

\smallskip

\noindent ({\it Point 3}) If $Z\sim\mathscr{N}(0,1)$, a simple application of the Fubini theorem (or, equivalently, an integration by parts) yields that $E[Zf(Z)]= E[f'(Z)]$ for every smooth $f$. Now suppose that $E[Zf(Z) - f'(Z)]=0$ for every function $f$ as in the statement, so that this equality holds in particular for $f=f_z$ and for every $z\in\R$. By integrating both sides of (\ref{stein-equation}) with respect to the law of $Z$, this yields that $P(Z\leq z) =P(N\leq z) $ for every $z\in\R$, and therefore that $Z$ and $N$ have the same law.
\fin

\begin{rem}{\rm
Formulae (\ref{stein-equation}) and (\ref{StCrudeBound}) are known, respectively, as {\sl Stein's equation} and {\sl Stein's bound}. As already evoked in the Introduction, Point 3 in the statement of Lemma \ref{stein} is customarily referred to as {\sl Stein's lemma}.
}
\end{rem}

\subsection{General bounds on the Kolmogorov distance}

We now face the problem of establishing a bound on the normal approximation of a centered and Malliavin-differentiable random variable. The next statement contains one of the main findings of \cite{NP-PTRF}.

\begin{thm}[See \cite{NP-PTRF}]\label{nou-pec}
Let $Z\in\mathbb{D}^{1,2}$ be such that $E(Z)=0$
and ${\rm Var}(Z)=1$. Then, for $N\sim \mathscr{N}(0,1)$,
\begin{equation}\label{GioIvan}
d_{Kol}(Z,N)
\leq \sqrt{{\rm Var}\big(\langle DZ,-DL^{-1}Z\rangle_\HH\big)}.
\end{equation}
\end{thm}
{\it Proof}.
In view of (\ref{StCrudeBound}), it is enough to prove that, for every Lipschitz function $f$ with Lipschitz constant less or equal to 1, one has that the quantity $|E[Zf(Z) - f'(Z)]|$ is less or equal to the RHS of (\ref{GioIvan}).
Start by considering a function $f:\R\to\R$ which is $C^1$ and such that $\|f'\|_\infty\leq 1$.
Relation (\ref{ivangioformula}) yields
\[
E\big(Zf(Z)\big)= E\big(f'(Z)\langle DZ,-DL^{-1}Z\rangle_\HH\big),
\]
so that
\[
\big|E\big(f'(Z)\big)-E\big(Zf(Z)\big)\big|=\big|E\big(f'(Z)(1-\langle DZ,-DL^{-1}Z\rangle_\HH)\big)\big|
\leq E\big|1-\langle DZ,-DL^{-1}Z\rangle_\HH\big|.
\]
By a standard approximation argument (e.g. by using a convolution with an approximation of the
identity), one sees that the inequality $\big|E\big(f'(Z)\big)-E\big(Zf(Z)\big)\big|
\leq E\big|1-\langle DZ,-DL^{-1}Z\rangle_\HH\big|$ continues to hold when $f$ is Lipschitz with constant less or equal to 1.
Hence, by combining the previous estimates with (\ref{StCrudeBound}), we infer that
\begin{eqnarray*}
d_{Kol} (Z,N)
\leq  E\big|1-\langle DZ,-DL^{-1}Z\rangle_\HH\big| \leq  \sqrt{E\big(1-\langle DZ,-DL^{-1}Z\rangle_\HH\big)^2}.
\end{eqnarray*}
Finally, the desired conclusion follows by observing that, if one chooses
$f(z)=z$ in (\ref{ivangioformula}), then one obtains
\begin{equation}\label{blabla}
E(\langle DZ,-DL^{-1}Z\rangle_\HH)=E(Z^2)=1,
\end{equation}
so that
$
E\left[\big(1-\langle DZ,-DL^{-1}Z\rangle_\HH\big)^2\right]={\rm Var}\big(\langle DZ,-DL^{-1}Z\rangle_\HH\big).
$
\fin

\begin{rem}{\rm By using the standard properties of conditional expectations, one sees that (\ref{GioIvan}) also implies the ``finer'' bound
\begin{equation}\label{GioIvan2}
d_{Kol}(Z,N)\leq \sqrt{{\rm Var}\big(g(Z)\big)},
\end{equation}
where $g(Z) = E[\langle DZ,-DL^{-1}Z\rangle_\HH | Z]$. In general, it is quite difficult to obtain an explicit expression of the function $g$. However, if some crude estimates on $g$ are available, then one can obtain explicit upper and lower bounds for the densities and the tail probabilities of the random variable $Z$. The reader is referred to Nourdin and Viens \cite{NV} and Viens \cite{Viens} for several results in this direction, and to Breton et al. \cite{breton-nourdin-peccati} for some statistical applications of these ideas.
}
\end{rem}

\subsection{Wiener chaos and the fourth moment condition}\label{SS : fourth}
In this section, we will apply Theorem \ref{nou-pec} to {\sl chaotic random variables}, that is, random variables having the special form of multiple Wiener-It\^{o} integrals of some fixed order $q\geq 2$. As announced in the Introduction, this allows to recover and refine some recent characterizations of CLTs on Wiener chaos (see \cite{NO, nunugio}).
We begin with a technical lemma.
\begin{lemma}\label{L : estimatesonvariance}
Fix an integer $q\geq 1$, and let $Z=I_q(f)$ (with $f\in\HH^{\odot q}$) be
such that ${\rm Var}(Z)=E(Z^2)=1$.
The following three identities are in order:
\begin{equation}\label{aa1}
\frac1q\|DZ\|^2_\HH-1=
q\sum_{r=1}^{q-1}(r-1)!\binom{q-1}{r-1}^2 I_{2q-2r}(f\widetilde{\otimes}_{r} f),
\end{equation}
\begin{equation}\label{lm-control-1d}
{\rm Var}\left(\frac1q\|DZ\|^2_\HH\right)=
\sum_{r=1}^{q-1}\frac{r^2}{q^2}\,r!^2\binom{q}{r}^4 (2q-2r)!\|f\widetilde{\otimes}_{r} f\|^2_{\HH^{\otimes 2q-2r}},
\end{equation}
and
\begin{equation}\label{aa3}
E(Z^4)-3=\frac{3}q\sum_{r=1}^{q-1}rr!^2\binom{q}{r}^4(2q-2r)!
\|f\widetilde{\otimes}_{r} f\|^2_{\HH^{\otimes 2q-2r}}.
\end{equation}
In particular,
\begin{equation}\label{momentdordre4}
{\rm Var}\left(\frac1q\|DZ\|^2_\HH\right)\leq \frac{q-1}{3q}\big(E(Z^4)-3\big).
\end{equation}
\end{lemma}
{\it Proof}.
Without loss of generality, we can assume that
$\HH$ is equal to $L^{2}(A,\mathcal{A},\nu)$, where $({A},\mathcal{A})$ is
a measurable space and $\nu$ a $\sigma$-finite measure  without atoms.
For any $a\in{A}$, we have
$D_aZ=qI_{q-1}\big(f(\cdot,a)\big)$ so that
\begin{eqnarray*}
\frac1q\|DZ\|^2_\HH&=&q\int_{{A}}
I_{q-1}\big(f(\cdot,a)\big)^2\nu(da)\\
&=&q\int_{A}\sum_{r=0}^{q-1}r!\binom{q-1}{r}^2 I_{2q-2-2r}\big(f(\cdot,a)\otimes_r f(\cdot,a)\big)\nu(da)
\quad
\mbox{by (\ref{multiplication})}\\
&=&q\sum_{r=0}^{q-1}r!\binom{q-1}{r}^2 I_{2q-2-2r}\left(\int_{A}f(\cdot,a)\otimes_r f(\cdot,a)\nu(da)\right)\notag\\
&=&q\sum_{r=0}^{q-1}r!\binom{q-1}{r}^2 I_{2q-2-2r}(f\otimes_{r+1} f)\\
&=&q\sum_{r=1}^{q}(r-1)!\binom{q-1}{r-1}^2 I_{2q-2r}(f\otimes_{r} f).\\
&=&q!\|f\|^2_{\HH^{\otimes q}}+q\sum_{r=1}^{q-1}(r-1)!\binom{q-1}{r-1}^2 I_{2q-2r}(f\otimes_{r} f).
\end{eqnarray*}
Since $E(Z^2)=q!\|f\|^2_{\HH^{\otimes q}}$, the proof of (\ref{aa1}) is finished.
The identity (\ref{lm-control-1d}) follows from (\ref{aa1}) and the orthogonality properties of multiple stochastic integrals.
Using (in order) formula (\ref{k1}) and the relation $D(Z^3)=3Z^2DZ$, we infer that
\begin{equation}
E(Z^4)=\frac1qE\big(\delta DZ\times Z^3\big)
=\frac1q E\big(\langle DZ,D(Z^3)\rangle_\HH\big)
=\frac3q E\big(Z^2\|DZ\|^2_\HH\big).\label{aa2}
\end{equation}
Moreover, the multiplication formula (\ref{multiplication}) yields
\begin{equation}
Z^2=I_q(f)^2=\sum_{s=0}^{q}s!\binom{q}{s}^2 I_{2q-2s}(f\otimes_s f).
\label{dollar2}
\end{equation}
By combining this last identity with (\ref{aa1}) and (\ref{aa2}), we obtain (\ref{aa3}) and finally (\ref{momentdordre4}).
\fin


As a consequence of Lemma \ref{L : estimatesonvariance}, we deduce the following bound on the Kolmogorov distance -- first proved in \cite{Noupecrei3}.
\begin{thm}[See \cite{Noupecrei3}] \label{mult-thm}
Let $Z$ belong to the $q$th chaos $\mathcal{H}_q$ of $X$, for some $q\geq 2$.
Suppose moreover that ${\rm Var}(Z)=E(Z^2)=1$.
Then
\begin{equation}\label{mult-ineq}
d_{Kol}(Z,N)\leq  \sqrt{\frac{q-1}{3q}\big(E(Z^4)-3\big)}.
\end{equation}
\end{thm}
{\it Proof}. Since
$L^{-1}Z=-\frac1qZ$, we have $\langle DZ,-DL^{-1}Z\rangle_\HH=\frac1q\|DZ\|^2_\HH$.
So, we only need to apply Theorem \ref{nou-pec} and formula (\ref{momentdordre4}).  \fin

The estimate (\ref{mult-ineq}) allows to deduce the following characterization of CLTs on Wiener chaos. Note that the equivalence of Point (i) and Point (ii) in the next statement was first proved by Nualart and Peccati in \cite{nunugio} (by completely different techniques based on stochastic time-changes), whereas the equivalence of Point (iii) was first obtained by Nualart and Ortiz-Latorre in \cite{NO} (by means of Malliavin calculus, but not of Stein's method).

\begin{thm}[see \cite{NO, nunugio}]\label{T : NPNOPTP}
Let $(Z_n)$ be a sequence of random variables belonging to the $q$th chaos $\mathcal{H}_q$
of $X$, for some fixed $q\geq 2$. Assume that
${\rm Var}(Z_n)=E(Z_n^2)=1$ for all $n$.
Then,
as $n\to\infty$, the following three assertions are equivalent:
\begin{enumerate}
\item[\rm (i)] $Z_n\overset{{\rm Law}}{\longrightarrow}N\sim\mathscr{N}(0,1)$;
\item[\rm (ii)] $E(Z_n^4)\to E(N^4) = 3$;
\item[\rm (iii)] ${\rm Var}\left(\frac1q\|DZ_n\|^2_\HH\right)\to 0$.
\end{enumerate}
\end{thm}
{\it Proof}.
For every $n$, write $Z_n=I_q(f_n)$ with $f_n\in\HH^{\odot q}$ uniquely determined.
The implication (iii) $\to$ (i) is a direct application of Theorem \ref{mult-thm}, and of the fact
that the topology of the Kolmogorov distance is stronger than the topology of the convergence
in law.
The implication (i) $\to$ (ii) comes from a bounded convergence argument (observe
that $\sup_{n\geq 1}E(Z_n^4)<\infty$ by the hypercontractivity relation (\ref{hypercontractivity})).
Finally, let us prove the implication (ii) $\to$ (iii).
Suppose that (ii) is in order. Then, by virtue of (\ref{aa3}), we have that
$\|f_n\widetilde{\otimes}_{r} f_n\|%
_{\HH^{\otimes 2q-2r}}$
tends to zero, as $n\to\infty$, for all (fixed) $r\in\{1,\ldots,q-1\}$.
Hence, (\ref{lm-control-1d}) allows to conclude that
(iii) is in order.
The proof of Theorem \ref{T : NPNOPTP} is thus complete.\fin

\begin{rem}{\rm
Theorem \ref{T : NPNOPTP} has been applied to a variety of situations: see e.g.
(but the list is by no means exhaustive)
Barndorff-Nielsen et al. \cite{Barndorff...},
Corcuera et al. \cite{CorNuWoe},
Marinucci and Peccati \cite{MaPeAb},
Neuenkirch and Nourdin \cite{NeuNourdin},
Nourdin and Peccati \cite{NouPecIBM}
and Tudor and Viens \cite{TudorViens}, and the references therein. See Peccati and Taqqu
\cite{PecTaq_SURV} for several combinatorial interpretations of these results.
}
\end{rem}

By combining Theorem \ref{mult-thm} and Theorem \ref{T : NPNOPTP},
we obtain the following result.

\begin{cor}
Let the assumptions of Corollary \ref{T : NPNOPTP} prevail. As $n\to\infty$, the following assertions are equivalent:
\begin{enumerate}
\item[\rm (a)] $Z_n\overset{{\rm Law}}{\longrightarrow}N\sim \mathscr{N}(0,1)$;
\item[\rm (b)] $d_{Kol}(Z_n,N)\to 0$.
\end{enumerate}
\end{cor}
{\it Proof}. Of course, only the implication (a) $\to$ (b) has to be proved. Assume that (a) is in order.
By Corollary \ref{T : NPNOPTP}, we have that ${\rm Var}\left(\frac1q\|DZ_n\|^2_\HH\right)\to 0$.
Using Theorem \ref{mult-thm}, we get that (b) holds, and the proof is done.\fin

\subsection{Quadratic variation of the fractional Brownian motion, part one}\label{BreuMaj}

In this section, we
use Theorem \ref{nou-pec} in order to derive an explicit bound for the
second-order approximation of the quadratic variation of a fractional Brownian motion.

Let $B=\{B_t:t\geq 0\}$ be a fractional Brownian motion with Hurst
index $H\in(0,1)$. This means that $B$ is a centered Gaussian process, started from zero and with covariance function $E%
(B_{s}B_{t})=R(s,t)$ given by
\[
R(s,t)=\frac{1}{2}\left( t^{2H}+s^{2H}-|t-s|^{2H}\right),\quad
s,t\geq 0.
\]
The fractional Brownian motion of index $H$ is the \textsl{only} centered Gaussian processes normalized in such a way that ${\rm Var}(B_1)=1$,
and such that $B$ is selfsimilar with index $H$ and has stationary increments.
If $H=1/2$ then $R(s,t)=\min(s,t)$ and $B$ is simply a standard Brownian motion.
If $H\neq 1/2$, then $B$ is neither a (semi)martingale nor a Markov process (see e.g.
\cite{nualartbook} for more details).

As already explained in the Introduction (see Example \ref{Exemple : cov}), for any choice of the Hurst parameter $H\in(0,1)$ the Gaussian space generated by $B$ can be identified
with an isonormal Gaussian process $X=\{X(h):h\in\HH\}$, where the real and separable Hilbert space
$\EuFrak H$ is defined as follows: (i) denote by $\mathscr{E}$ the
set of all $\mathbb{R}$-valued step functions on $[0,\infty)$, (ii)
define $\EuFrak H$ as the Hilbert space obtained by closing
$\mathscr{E}$ with respect to the scalar product
\[
\left\langle
{\mathbf{1}}_{[0,t]},{\mathbf{1}}_{[0,s]}\right\rangle _{\EuFrak
H}=R(t,s).
\]
In particular, with such a notation, one has that
$B_t=X(\mathbf{1}_{[0,t]})$.

Set
\[
Z_n=\frac{1}{\sigma_n}\sum_{k=0}^{n -1}\big[(B_{k+1}-B_k)^2-1\big]
\overset{\rm Law}{=}\frac{n^{2H}}{\sigma_n}\sum_{k=0}^{n -1}\big[(B_{(k+1)/n}-B_{k/n})^2-n^{-2H}\big]
\]
where $\sigma_n>0$ is chosen so that $E(Z_n^2)=1$. It is well-known (see e.g. \cite{BrMa}) that, for every $H\leq 3/4$ and for $n\to\infty$, one has that $Z_n$ converges in law to $N\sim\mathscr{N}(0,1)$. The following result uses Stein's method in order to obtain an explicit bound for the
Kolmogorov distance between $Z_n$ and $N$. It was first proved in \cite{NP-PTRF} (for the case $H<3/4$) and \cite{breton-nourdin} (for $H=3/4$).

\begin{thm}\label{BM}
Let $N\sim \mathscr{N}(0,1)$ and assume that $H\leq 3/4$. Then, there exists a
constant $c_H>0$ (depending only on $H$) such that, for every
$n\geq 1$,
\begin{equation}\label{bmstat}
d_{Kol}(Z_n,N)\leq c_H\times
\left\{\begin{array}{lll}
\frac1{\sqrt{n}}&\,\,\mbox{if $H\in (0,\frac12]$}\\
\\
n^{2H-\frac32}&\,\,\mbox{if $H\in [\frac12,\frac34)$}\\
\\
\frac1{\sqrt{\log n}}&\,\,\mbox{if $H=\frac34$}\\
\end{array}\right..
\end{equation}
\end{thm}

\begin{remark}
{\rm
\begin{enumerate}
\item By inspection of the forthcoming proof of Theorem \ref{BM}, one sees that
$\lim_{n\to\infty}\frac{\sigma_n^2}{n}=2\sum_{r\in\mathbb{Z}}\rho^2(r)$ if $H\in(0,3/4)$,
with $\rho$ given by (\ref{rho}), and
$\lim_{n\to\infty}\frac{\sigma_n^2}{n\log n}=9/16$ if $H=3/4$.
\item When $H>3/4$, the sequence $(Z_n)$ does not converge in law to $\mathscr{N}(0,1)$.
Actually, $Z_n \underset{n\to\infty}{\overset{{\rm Law}}{\longrightarrow}} Z_\infty\sim\mbox{``Hermite random variable''}$
and, using a result by Davydov and Martynova \cite{davydov-martynova}, one can also associate a bound to this convergence.
See \cite{breton-nourdin} for details on this result.
\item More generally, and using the analogous computations, one can associate bounds with the convergence of sequence
\[
Z^{(q)}_n=\frac{1}{\sigma^{(q)}_n}\sum_{k=0}^{n-1}H_q(B_{k+1}-B_k)
\overset{\rm Law}{=}\frac{1}{\sigma^{(q)}_n}\sum_{k=0}^{n -1}H_q(n^H(B_{(k+1)/n}-B_{k/n})\big)
\]
towards $N\sim\mathscr{N}(0,1)$, where
$H_q$ ($q\geq 3$) denotes the $q$th Hermite polynomial (as defined in (\ref{hermite})), and $\sigma^{(q)}_n$ is some appropriate normalizing constant. In this case, the critical value is $H=1-1/(2q)$ instead of $H=3/4$. See \cite{NP-PTRF} for details.
\end{enumerate}
}
\end{remark}

In order to show Theorem \ref{BM}, we will need the following ancillary result,
whose proof is obvious and left to the reader.
\begin{lemma}\label{lm-rho}
\begin{itemize}
\item[\rm 1.] For $r\in\mathbb{Z}$, let
\begin{equation}\label{rho}
\rho(r)=\frac12\big(|r+1|^{2H}+|r-1|^{2H}-2|r|^{2H}\big).
\end{equation}
If $H\neq \frac12$, one has $\rho(r)\sim H(2H-1)|r|^{2H-2}$ as $|r|\to\infty$.
If $H=\frac12$ and $|r|\geq 1$, one has $\rho(r)=0$.
Consequently, $\sum_{r\in\mathbb{Z}} \rho^2(r)<\infty$ if and only if $H<3/4$.
\item[\rm 2.] For all $\alpha>-1$,
we have $\sum_{r=1}^{n-1} r^\alpha \sim n^{\alpha+1}/(\alpha+1)$ as $n\to\infty$.
\end{itemize}
\end{lemma}

\medskip

We are now ready to prove the main result of this section.

\bigskip

\noindent{\it Proof of Theorem \ref{BM}}.
Since $\|
{\bf 1}_{[k,k+1]}
\|^2_\HH=E\big((B_{k+1}-B_k)^2\big)=1$, we have,
by (\ref{hermite-link}),
\[
(B_{k+1}-B_k)^2-1=I_2({\bf 1}_{[k,k+1]}^{\otimes 2})
\]
so that $Z_n=I_2(f_n)$ with $f_n=\frac{1}{\sigma_n}\sum_{k=0}^{n-1}{\bf 1}_{[k,k+1]}^{\otimes 2}\in\HH^{\odot 2}$.
Let us compute the exact value of $\sigma_n$.
Observe that $\langle {\bf 1}_{[k,k+1]},{\bf 1}_{[l,l+1]}\rangle_\HH
=E\big((B_{k+1}-B_k)(B_{l+1}-B_l)\big)=\rho(k-l)$ with $\rho$ given by (\ref{rho}).
Hence
\begin{eqnarray*}
&&E\left[\left(\sum_{k=0}^{n -1}\big[(B_{k+1}-B_k)^2-1\big]\right)^2\right]\\
&=&E\left[\left(\sum_{k=0}^{n -1} I_2({\bf 1}_{[k,k+1]}^{\otimes 2})\right)^2\right]
=\sum_{k,l=0}^{n-1}E\big[I_2({\bf 1}_{[k,k+1]}^{\otimes 2})I_2({\bf 1}_{[l,l+1]}^{\otimes 2})\big]\\
&=&2\sum_{k,l=0}^{n-1}\langle {\bf 1}_{[k,k+1]},{\bf 1}_{[l,l+1]}\rangle_\HH^2=2\sum_{k,l=0}^{n-1}\rho^2(k-l).
\end{eqnarray*}
That is,
\begin{eqnarray*}
\sigma_n ^2= 2\sum_{k,l=0}^{n-1}\rho^2(k-l)
=2\sum_{l=0}^{n-1}\sum_{r=-l}^{n-1-l}\rho^2(r)
=2\left(n\sum_{|r|<n} \rho^2(r) - \sum_{|r|<n} \big(|r|+1\big)
\rho^2(r)\right).
\end{eqnarray*}
Assume that $H<3/4$. Then, we have
\[
\frac{\sigma_n^2}{n}=2\sum_{r\in\mathbb{Z}} \rho^2(r)\left(1-\frac{|r|+1}{n}\right)
{\bf 1}_{\{|r|<n\}}.
\]
Since $\sum_{r\in\mathbb{Z}} \rho^2(r)<\infty$, we obtain, by bounded Lebesgue convergence:
\begin{equation}\label{sigman}
\lim_{n\to\infty}
\frac{\sigma_n^2}{n}=2\sum_{r\in\mathbb{Z}} \rho^2(r).
\end{equation}
Assume that $H=3/4$. We have
$
\rho^2(r)\sim \frac{9}{64|r|}
$ as $|r|\to\infty$.
Therefore, as $n\to\infty$,
\[
n\sum_{|r|<n} \rho^2(r) \sim
\frac{9n}{64}\sum_{0<|r|<n}\frac{1}{|r|}
\sim \frac{9n\log n}{32}
\]
and
\[
\sum_{|r|<n} \big(|r|+1\big)\rho^2(r) \sim
\frac{9}{64}\sum_{|r|<n} 1
\sim \frac{9n}{32}.
\]
We deduce that
\begin{equation}\label{sigman2}
\lim_{n\to\infty}\frac{\sigma_n^2}{n\log n}=\frac{9}{16}.
\end{equation}
Now, we have, see (\ref{lm-control-1d}) for the first equality,
\begin{eqnarray*}
{\rm Var}\big(\frac12\|DZ_n\|^2_\HH\big)&=&
\frac12\|f_n\otimes_1 f_n\|^2_{\HH^{\otimes 2}}=\frac1{2\sigma_n^4}\left\| \sum_{k,l=0}^{n-1}
{\bf 1}_{[k,k+1]}^{\otimes 2}\otimes_1{\bf 1}_{[l,l+1]}^{\otimes 2}\right\|_\HH^2\\
&=&\frac1{2\sigma_n^4}\left\| \sum_{k,l=0}^{n-1}\rho(k-l)
{\bf 1}_{[k,k+1]}\otimes{\bf 1}_{[l,l+1]}\right\|_\HH^2\\
&=&\frac1{2\sigma_n^4}\sum_{i,j,k,l=0}^{n-1}\rho(k-l)\rho(i-j)\rho(k-i)\rho(l-j)\\
&\leq&\frac1{4\sigma_n^4}\sum_{i,j,k,l=0}^{n-1}|\rho(k-i)||\rho(i-j)|
\big(\rho^2(k-l)+\rho^2(l-j)\big)\\
&\leq&\frac1{2\sigma_n^4}\sum_{i,j,k=0}^{n-1}|\rho(k-i)||\rho(i-j)|\sum_{r=-n+1}^{n-1}
\rho^2(r)\\
&\leq&
\frac{n}{2\sigma_n^4}\left(
\sum_{s=-n+1}^{n-1}
|\rho(s)|\right)^2 \sum_{r=-n+1}^{n-1}
\rho^2(r).
\end{eqnarray*}
If $H\leq 1/2$ then $\sum_{s\in\mathbb{Z}}|\rho(s)|<\infty$
and $\sum_{r\in\mathbb{Z}}\rho^2(r)<\infty$
so that, in view of (\ref{sigman}),
${\rm Var}\big(\frac12\|DZ_n\|^2_\HH\big)=O(n^{-1})$.
If $1/2<H< 3/4$ then $\sum_{s=-n+1}^{n-1}|\rho(s)|=O(n^{2H-1})$ (see Lemma \ref{lm-rho})
and $\sum_{r\in\mathbb{Z}}\rho^2(r)<\infty$
so that, in view of (\ref{sigman}),
 one has ${\rm Var}\big(\frac12\|DZ_n\|^2_\HH\big)=O(n^{4H-3})$.
If $H=3/4$ then $\sum_{s=-n+1}^{n-1}|\rho(s)|=O(\sqrt{n})$
and $\sum_{r=-n+1}^{n-1}\rho^2(r)=O(\log n)$ (indeed, by Lemma \ref{lm-rho},
$\rho^2(r)\sim \frac{{\rm cst}}{|r|}$
as $|r|\to\infty$)
so that, in view of (\ref{sigman2}),
${\rm Var}\big(\frac12\|DZ_n\|^2_\HH\big)=O(1/\log n)$.
Finally, the desired conclusion follows from Theorem \ref{mult-thm}.
\fin

\subsection{The method of (fourth) moments: explicit estimates via interpolation}
It is clear that the combination of Theorem \ref{mult-thm} and Theorem \ref{T : NPNOPTP} provides a remarkable simplification of the method of moments and cumulants, as applied to the derivation of CLTs on a fixed Wiener chaos (further generalizations of these results, concerning in particular multi-dimensional CLTs, are discussed in the forthcoming Section \ref{S : MULTI}). In particular, one deduces from (\ref{mult-ineq}) that, for a sequence of chaotic random variables with unit variance, the speed of convergence to zero of the fourth cumulants $E(Z_n^4)-3$ also determines the speed of convergence in the Kolmogorov distance.

In this section, we shall state and prove a new upper bound, showing that, for a normalized chaotic sequence $\{Z_n : n\geq 1\}$ converging in distribution to $N\sim\mathscr{N}(0,1)$, the convergence to zero of $E(Z_n^k) - E(N^k)$ is always dominated by the speed of convergence of the square root of $E(Z_n^4) - E(N^4) = E(Z_n^4) - 3$. To do this, we shall apply a well-known Gaussian interpolation technique, which has been essentially introduced by Talagrand (see e.g. \cite{tala}); note that a similar approach has recently been adopted in \cite{Noupecrei3}, in order to deduce a universal characterization of CLTs for sequences of homogeneous sums.

\begin{rem}{\rm
\begin{itemize}
\item[1.] In principle, 
one could deduce from the results of this section that, for every $k\geq 3$, the speed of convergence to zero of $k$th cumulant of $Z_n$ is always dominated by the speed of convergence of the fourth cumulant $E(Z_n^4)-3$.
\item[2.] We recall that the explicit computation of moments and cumulants of chaotic random variables is often performed by means of a class of combinatorial devices, known as {\sl diagram formulae}. This tools are not needed in our analysis, as we rather rely on multiplication formulae and integration by parts techniques from Malliavin calculus. See \cite[Section 3]{PecTaq_SURV} for a recent and self-contained introduction to moments, cumulants and diagram formulae.
\end{itemize}
}
\end{rem}

\begin{prop}\label{P : moments}
Let $q\geq 2$ be an integer, and let $Z$ be an element of the $q$th chaos $\mathcal{H}_q$
of $X$. Assume that ${\rm Var}(Z)=E(Z^2)=1$, and
let $N\sim \mathscr{N}(0,1)$.
Then, for all integer $k\geq 3$,
\begin{equation}\label{moment}
\big|E(Z^k) -E(N^k)\big| \leq
c_{k,q}\sqrt{E(Z^4)-E(N^4)},
\end{equation}
where the constant $c_{k,q}$ is given by
\[
c_{k,q}=(k-1)2^{k-\frac52}
\sqrt{\frac{q-1}{3q}}
\left(
\sqrt{\frac{(2k-4)!}{2^{k-2}(k-2)!}}
+(2k-5)^{\frac{kq}2-q}
\right).
\]
\end{prop}
{\it Proof}. Without loss of generality, we can assume that $N$ is independent of the underlying
isonormal Gaussian process
$X$. Fix an integer $k\geq 3$.
By denoting $\Psi(t)=E\big[(\sqrt{1-t}Z+\sqrt{t}N)^k\big]$, $t\in[0,1]$, we have
\[
\big|E(Z^k) -E(N^k)\big|=\big|\Psi(1)
-\Psi(0)\big|\leq \int_0^1 |\Psi'(t)|dt,
\]
where the derivative $\Psi'$ is easily seen to exist on $(0,1)$, and moreover one has
\[
\Psi'(t)=\frac{k}{2\sqrt{t}} E\big[(\sqrt{1-t}Z+\sqrt{t}N)^{k-1}N\big]-
\frac{k}{2\sqrt{1-t}}E\big[(\sqrt{1-t}Z+\sqrt{t}N)^{k-1}Z\big].
\]
By integrating by parts and by using the explicit expression of the Gaussian density, one infers that
\begin{eqnarray*}
E\big[(\sqrt{1-t}Z+\sqrt{t}N)^{k-1}N\big]&=&E\left[E\big[(\sqrt{1-t}z+\sqrt{t}N)^{k-1}N\big]_{|z=Z}\right]\\
&=&(k-1)\sqrt{t}\,E\left[E\big[(\sqrt{1-t}z+\sqrt{t}N)^{k-2}\big]_{|z=Z}\right]\\
&=&(k-1)\sqrt{t}\,E\big[(\sqrt{1-t}Z+\sqrt{t}N)^{k-2}\big].
\end{eqnarray*}
Similarly, using this time (\ref{ivangioformula}) in order to perform the integration by parts and taking into account that
$\langle DZ,-DL^{-1}Z\rangle_\HH=\frac1q\|DZ\|^2_{\HH}$ because $Z\in\mathcal{H}_q$, we can write
\begin{eqnarray*}
E\big[(\sqrt{1-t}Z+\sqrt{t}N)^{k-1}Z\big]&=&E\left[E\big[(\sqrt{1-t}Z+\sqrt{t}x)^{k-1}Z\big]_{|x=N}\right]\\
&=&(k-1)\sqrt{1-t}\,E\left[E\big[(\sqrt{1-t}Z+\sqrt{t}x)^{k-2}\frac1q\|DZ\|^2_{\HH}\big]_{|x=N}\right]\\
&=&(k-1)\sqrt{1-t}\,E\left[(\sqrt{1-t}Z+\sqrt{t}N)^{k-2}\frac1q\|DZ\|^2_{\HH}\right].
\end{eqnarray*}
Hence,
\[
\Psi'(t)=
\frac{k(k-1)}{2} E\left[\left(1-\frac1q\|DZ\|^2_{\HH}\right)(\sqrt{1-t}Z+\sqrt{t}N)^{k-2}\right],
\]
and consequently
\[
\big|\Psi'(t)\big|\leq
\frac{k(k-1)}{2} \sqrt{E\left[(\sqrt{1-t}Z+\sqrt{t}N)^{2k-4}\right]}\times
\sqrt{E\left[\left(1-\frac1q\|DZ\|^2_{\HH}\right)^2\right]}.
\]
By (\ref{blabla}) and (\ref{momentdordre4}), we have
\[
E\left[\left(1-\frac1q\|DZ\|^2_{\HH}\right)^2\right]
=
{\rm Var}\left(\frac1q\|DZ\|^2_{\HH}\right)\leq  \frac{q-1}{3q}\big(E(Z^4)-3\big).
\]
Using succesively $(x+y)^{2k-4}\leq 2^{2k-5}(x^{2k-4}+y^{2k-4})$,
$\sqrt{x+y}\leq\sqrt{x}+\sqrt{y}$,
inequality (\ref{hypercontractivity}) and $E(N^{2k-4})=(2k-4)!/(2^{k-2}(k-2)!)$,
we can write
\begin{eqnarray*}
\sqrt{E\left[(\sqrt{1-t}Z+\sqrt{t}N)^{2k-4}\right]}&\leq&2^{k-\frac52}(1-t)^{\frac{k}2-1}\sqrt{E(Z^{2k-4})}+2^{k-\frac52}t^{\frac{k}2-1}\sqrt{E(N^{2k-4})}\\
&\leq&2^{k-\frac52}(1-t)^{\frac{k}2-1}(2k-5)^{\frac{kq}2-q}\!+ \! 2^{k-\frac52}t^{\frac{k}2-1}\sqrt{\frac{(2k-4)!}{2^{k-2}(k-2)!}}\\
\end{eqnarray*}
so that
\[
\int_0^1 \sqrt{E\left[(\sqrt{1-t}Z+\sqrt{t}N)^{2k-4}\right]} dt \leq
\frac{2^{k-\frac32}}{k}\left[(2k-5)^{\frac{kq}2-q}+\sqrt{\frac{(2k-4)!}{2^{k-2}(k-2)!}}\right].
\]
Putting all these bounds together, one deduces the desired conclusion.
\fin

\section{Multidimensional case}\label{S : MULTI}

Here and for the rest of the section, we consider as given an isonormal Gaussian process $\{X(h) : h\in \EuFrak{H}\}$, over some real separable Hilbert space $\EuFrak{H}$.

\subsection{Main bounds}

We shall now present (without proof) a result taken from \cite{multivariate}, concerning the Gaussian approximation of vectors of random variables that are differentiable in the Malliavin sense. We recall that the {\sl Wasserstein distance} between the laws of two $\real^d$-valued
random vectors $X$ and $Y$, noted $d_{\rm W}(X,Y)$, is given by
$$ d_{\rm W}(X,Y):=\sup_{g\in\mathscr{H}; \|g\|_{Lip}\leq 1} \big\vert E[g(X)]-E[g(Y)] \big\vert,$$
where $\mathscr{H}$ indicates the class of all Lipschitz functions,
that is, the collection of all functions $g:\real^d\to\real$ such
that $$\displaystyle{\|g\|_{Lip}:=\sup_{x\neq y}\frac{\vert
g(x)-g(y) \vert}{\| x-y \|_{\real^d}}<\infty}$$
(with $\|\cdot\|_{\real^d}$ the usual Euclidian norm on $\real^d$). Also, we recall that the {\sl operator norm}
of a $ d\times d$ matrix $A$ over $\real$ is given by
$ \|A\|_{op} :=\sup_{\|x\|_{\real^d}=1}\|A
x\|_{\real^d}.$

Note that, in the following statement, we require that the approximating Gaussian vector has a positive definite covariance matrix.

\begin{theorem}[See \cite{multivariate}]
\label{theo:majDist} Fix $d\geq 2$ and let $C=(C_{ij})_{1\leq i,j\leq d}$ be a $d\times d$ positive definite matrix. Suppose
that $N\sim\mathscr{N}_d(0,C)$, 
and assume that $Z=(Z_1,\ldots,Z_d)$ is a
$\real^d$-valued random vector such that $E[Z_i]=0$ and $Z_i \in
\mathbb{D}^{1,2}$ for every $i=1,\ldots,d$. Then,
\begin{eqnarray}
 d_{\rm W}(Z,N)
&\leq& \| C^{-1}\|_{op} \,\, \| C\|_{op}^{1/2} \sqrt{\sum_{i,j=1}^d
E[( C_{ij} - \langle DZ_i, -DL^{-1}Z_j \rangle_{\EuFrak H})^2] }.
\notag
\end{eqnarray}
\end{theorem}

In what follows, we shall use once again interpolation techniques in order to
partially generalize Theorem \ref{theo:majDist} to the case where the approximating
covariance matrix $C$ is not necessarily positive definite. This additional difficulty forces us to work with functions that are smoother than the ones involved in the definition of the Wasserstein distance. To this end, we will adopt the following simplified notation: for every $\varphi:\R^d\to\R$ of class $C^2$, we set
\begin{eqnarray*}
\|\varphi''\|_\infty=\max_{i,j=1,\ldots,d}
\sup_{z\in\R^d}\left|
\frac{\partial^2\varphi}{\partial x_i\partial x_j}(z)\right|.
\end{eqnarray*}

\begin{theorem}[See \cite{Noupecrei3}] \label{theo:majDist2}
Fix $d\geq 2$, and let $C=(C_{ij})_{1\leq i,j\leq d}$ be a $d\times d$ covariance matrix. Suppose that
$N\sim \mathscr{N}_d(0,C)$ and that $Z=(Z_1,\ldots, Z_d)$ is a $\R^d$-valued random vector such that $E[Z_i]=0$ and $Z_i \in
\mathbb{D}^{1,2}$ for every $i=1,\ldots,d$. Then, for every $\varphi:\R^d\to \R$ belonging to $C^2$ such that $\|\varphi''\|_\infty<\infty$, we have
\begin{equation}\label{WassTWO}
\big|E[\varphi(Z)]-E[\varphi(N)]\big|\leq
\frac12\|\varphi''\|_\infty
\sum_{i,j=1}^d
E\left[\left|C_{i,j}-\langle DZ_j,-DL^{-1}Z_i\rangle_\HH\right|\right].
\end{equation}
\end{theorem}
{\it Proof}.
Without loss of generality, we assume that $N$ is independent of
the underlying isonormal Gaussian process $X$.
Let $\varphi:\R^d\to \R$ be a $C^2$-function such that $\|\varphi''\|_\infty<\infty$.
For any $t\in[0,1]$, set
$
\Psi(t)=E\big[\varphi\big(\sqrt{1-t}Z+\sqrt{t}N\big)\big],
$
so that
\[
\big|E[\varphi(Z)]-E[\varphi(N)]\big|=\big|\Psi(1)
-\Psi(0)\big|\leq \int_0^1 |\Psi'(t)|dt.
\]
We easily see that $\Psi$ is differentiable on $(0,1)$ with
\[
\Psi'(t)=
\sum_{i=1}^d E\left[\frac{\partial\varphi}{\partial x_i}\big(\sqrt{1-t}Z
+\sqrt{t}N\big)\left(
\frac1{2\sqrt{t}}N_i-\frac{1}{2\sqrt{1-t}}Z_i\right)
\right].
\]
By integrating by parts, we can write
\begin{eqnarray*}
&&E\left[
\frac{\partial\varphi}{\partial x_i}\big(\sqrt{1-t}Z
+\sqrt{t}N\big)N_i
\right]\\
&=&E\left\{
E\left[
\frac{\partial\varphi}{\partial x_i}\big(\sqrt{1-t}z
+\sqrt{t}N\big)N_i
\right]_{|z=Z}\right\}\\
&=&\sqrt{t}\sum_{j=1}^d C_{i,j}\,E\left\{
E\left[
\frac{\partial^2\varphi}{\partial x_i\partial x_j}\big(\sqrt{1-t}z
+\sqrt{t}N\big)
\right]_{|z=Z}\right\}\\
&=&\sqrt{t}\sum_{j=1}^d C_{i,j}\,
E\left[
\frac{\partial^2\varphi}{\partial x_i\partial x_j}\big(\sqrt{1-t}Z
+\sqrt{t}N\big)
\right].
\end{eqnarray*}
By using (\ref{ivangioformula}) in order to perform the integration by parts, we can also write
\begin{eqnarray*}
&&E\left[
\frac{\partial\varphi}{\partial x_i}\big(\sqrt{1-t}Z
+\sqrt{t}N\big)Z_i
\right]\\
&=&E\left\{
E\left[
\frac{\partial\varphi}{\partial x_i}\big(\sqrt{1-t}Z
+\sqrt{t}x\big)Z_i
\right]_{|x=N}\right\}\\
&=&\sqrt{1-t}\sum_{j=1}^d E\left\{
E\left[
\frac{\partial^2\varphi}{\partial x_i\partial x_j}\big(\sqrt{1-t}Z
+\sqrt{t}x\big)
\langle DZ_j,-DL^{-1}Z_i\rangle_\HH
\right]_{|x=N}\right\}\\
&=&\sqrt{1-t}\sum_{j=1}^d
E\left[
\frac{\partial^2\varphi}{\partial x_i\partial x_j}\big(\sqrt{1-t}Z
+\sqrt{t}N\big)
\langle DZ_j,-DL^{-1}Z_i\rangle_\HH
\right].
\end{eqnarray*}
Hence
\[
\Psi'(t)=\frac1{2}
\sum_{i,j=1}^d
E\left[
\frac{\partial^2\varphi}{\partial x_i\partial x_j}\big(\sqrt{1-t}Z
+\sqrt{t}N\big)\left(C_{i,j}-
\langle DZ_j,-DL^{-1}Z_j\rangle_\HH\right)
\right],
\]
so that
\[
\int_0^1 |\Psi'(t)|dt\leq \frac12\|\varphi''\|_\infty
\sum_{i,j=1}^d
E\left[\left|C_{i,j}-\langle DZ_j,-DL^{-1}Z_i\rangle_\HH\right|\right]
\]
and the desired conclusion follows.
\fin

We now aim at applying Theorem
\ref{theo:majDist2} to vectors of multiple stochastic integrals.

\begin{corollary}
\label{prop:Chaos} Fix integers $d\geq 2$ and $1\leq q_1\leq\ldots\leq
q_d$. Consider a vector
$Z=(Z_1,\ldots,Z_d):=(I_{q_1}(f_1),\ldots,I_{q_d}(f_d))$ with
$f_{i}\in \EuFrak{H}^{\odot q_i}$ for any $i=1\ldots,d$. Let
$N\sim \mathcal{N}_d(0,C)$,
with
$C=(C_{ij})_{1\leq i,j\leq d}$ a $d\times d$ covariance matrix.
Then, for every $\varphi:\R^d\to \R$ belonging to $C^2$ such that $\|\varphi''\|_\infty<\infty$,
we have
\begin{equation}\label{star}
\big|E[\varphi(Z)]-E[\varphi(N)]\big|\leq
\frac12\|\varphi''\|_\infty
\sum_{i,j=1}^d
E\left[\left|C_{i,j}-\frac1{d_i}\langle DZ_j,DZ_i\rangle_\HH\right|\right].
\end{equation}
\end{corollary}
{\it Proof}.
We have $-L^{-1}Z_i=\frac{1}{d_i}\,Z_i$ so that
the desired conclusion follows from (\ref{WassTWO}).
\fin

When one applies Corollary \ref{prop:Chaos} in concrete situations, one can use
the following result in order to evaluate the right-hand side of (\ref{star}).
\begin{prop}\label{lm-control}
Let $F=I_p(f)$ and $G=I_q(g)$, with $f\in\HH^{\odot p}$ and
$g\in\HH^{\odot q}$ ($p,q\geq 1$). Let $a$ be a real constant. If
$p=q$, one has the estimate:
\begin{eqnarray*}
&&E\left[\left(a-\frac1p\left\langle DF,DG\right\rangle_\HH\right)^2\right] \leq (a-p!\langle f,g\rangle_{\HH^{\otimes p}})^2\\
&&\hskip2cm+
\frac{p^2}{2}\sum_{r=1}^{p-1}(r-1)!^2\binom{p-1}{r-1}^4(2p-2r)!
\big( \|f\otimes_{p-r}f\|^2_{\HH^{\otimes
2r}}+\|g\otimes_{p-r}g\|^2_{\HH^{\otimes 2r}}\big).
\end{eqnarray*}
On the other hand, if $p< q$, one has that
\begin{eqnarray*}
&&E\left[\left(a-\frac1q\left\langle DF,DG\right\rangle_\HH
\right)^2\right]
\leq a^2+p!^2\binom{q-1}{p-1}^2(q-p)!\|f\|^2_{\HH^{\otimes p}}\|g\otimes_{q-p}g\|_{\HH^{\otimes 2p}}\\
&&+
\frac{p^2}{2}\sum_{r=1}^{p-1}(r-1)!^2\binom{p-1}{r-1}^2\binom{q-1}{r-1}^2(p+q-2r)!\big(
\|f\otimes_{p-r}f\|^2_{\HH^{\otimes
2r}}+\|g\otimes_{q-r}g\|^2_{\HH^{\otimes 2r}}\big).
\end{eqnarray*}
\end{prop}

\begin{remark}\label{r}
{\rm
When bounding the right-hand side of (\ref{star}), we see that it
is sufficient to asses the quantities $\|f_i\otimes_r
f_i\|_{\HH^{\otimes2(q_i-r)}}$ for all $i=1,\ldots,d$ and
$r=1,\ldots,q_i-1$ on the one hand, and $E(Z_iZ_j)=q_i!\langle
f_i,f_j\rangle_{\HH^{\otimes q_i}}$ for all $i,j=1,\ldots,d$ such
that $q_i=q_j$ on the other hand. In particular, this fact allows to recover a result first proved by Peccati and Tudor in \cite{PTu04}, namely that, for vectors of multiple stochastic integrals whose covariance matrix is converging, \textsl{the componentwise convergence to a Gaussian distribution always implies joint convergence}.
}
\end{remark}

\noindent{\it Proof of Proposition \ref{lm-control}}. Without loss of
generality, we can assume that $\HH=L^{2}(A,\mathscr{A},\mu)$,
where $(A,\mathscr{A})$ is a measurable space, and $\mu$ is a
$\sigma$-finite and non-atomic measure. Thus, we can write
\begin{eqnarray*}
\langle DF,DG\rangle_\HH &=&p\,q\left\langle I_{p-1}(f),I_{q-1}(g)\right\rangle_\HH
=p\,q\int_A I_{p-1}\big(f(\cdot,t)\big)I_{q-1}\big(g(\cdot,t)\big)\mu(dt)\\
&=&p\,q\int_A \sum_{r=0}^{p\wedge q-1} r!\binom{p-1}{r} \binom{q-1}{r} I_{p+q-2-2r}\big(f(\cdot,t)\widetilde{\otimes}_r g(\cdot,t)\big)\mu(dt)\\
&=&p\,q\sum_{r=0}^{p\wedge q-1} r!\binom{p-1}{r}\binom{q-1}{r} I_{p+q-2-2r}(f\widetilde{\otimes}_{r+1}g)\\
&=&p\,q \sum_{r=1}^{p\wedge q} (r-1)!\binom{p-1}{r-1}\binom{q-1}{r-1} I_{p+q-2r}(f\widetilde{\otimes}_r g).
\end{eqnarray*}
It follows that
\begin{eqnarray}
&&E\left[\left(a-\frac1q\left\langle DF,DG\right\rangle_\HH\right)^2\right] \label{Murray}\\
&=&\left\lbrace
\begin{array}{l}
a^2+p^2\sum_{r=1}^{p}(r-1)!^2
\binom{p-1}{r-1}^2\binom{q-1}{r-1}^2 (p+q-2r)!
\|f\widetilde{\otimes}_r g\|^2_{\HH^{\otimes (p+q-2r)}} \textrm{ if } p< q,\\\\
(a-p!\langle f, g\rangle_{\EuFrak{H}^{\otimes
p}})^2+p^2\sum_{r=1}^{p-1}(r-1)!^2 \binom{p-1}{r-1}^4 (2p-2r)!
\|f\widetilde{\otimes}_r g\|^2_{\HH^{\otimes (2p-2r)}} \textrm{ if
} p=q.
\end{array}\notag
\right.
\end{eqnarray}
If $r<p\leq q$ then
\begin{eqnarray*}
\|f\widetilde{\otimes}_r g\|^2_{\HH^{\otimes (p+q-2r)}}
&\leq& \|f\otimes_r g\|^2_{\HH^{\otimes (p+q-2r)}}
=\langle f\otimes_{p-r} f, g\otimes_{q-r}g\rangle_{\HH^{\otimes 2r}}\\
&\leq&
\|f\otimes_{p-r}f\|_{\HH^{\otimes 2r}}\|g\otimes_{q-r}g\|_{\HH^{\otimes 2r}}\\
&\leq&\frac12\left(
\|f\otimes_{p-r}f\|_{\HH^{\otimes 2r}}^2+\|g\otimes_{q-r}g\|_{\HH^{\otimes 2r}}^2
\right).
\end{eqnarray*}
If $r=p<q$, then
$$
\|f\widetilde{\otimes}_p\, g\|^2_{\HH^{\otimes (q-p)}} \leq
\|f\otimes_p \,g\|^2_{\HH^{\otimes (q-p)}} \leq
\|f\|^2_{\HH^{\otimes p}}\|g\otimes_{q-p}g\|_{\HH^{\otimes 2p}}.
$$
If $r=p=q$, then $ f\widetilde{\otimes}_p g =\langle
f,g\rangle_{\HH^{\otimes p}}. $
By plugging these last expressions into (\ref{Murray}), we deduce
immediately the desired conclusion.
\fin

\subsection{Quadratic variation of fractional Brownian motion, continued}
In this section, we continue the example of Section \ref{BreuMaj}.
We still denote by $B$ a fractional Brownian
motion with Hurst index $H\in(0,3/4]$.
We set
\[
Z_n(t)=\frac{1}{\sigma_n}\sum_{k=0}^{\lfloor nt\rfloor -1}\big[(B_{k+1}-B_k)^2-1\big],\quad t\geq 0,
\]
where $\sigma_n>0$ is such that $E\big(Z_n(1)^2\big)=1$.
The following
statement contains the multidimensional
counterpart of Theorem \ref{BM}, namely a bound associated with the convergence of the finite dimensional distributions of $\{Z_n(t):\,t\geq 0\}$
towards a standard Brownian motion. A similar result can be of course recovered from Theorem
\ref{theo:majDist} -- see again \cite{multivariate}.
\begin{thm}\label{BM-rev}
Fix $d\geq 1$, and consider $0=t_0<t_1<\ldots<t_d$. Let $N\sim \mathscr{N}_d(0,I_d)$. There exists a
constant $c$ (depending only on $d$, $H$ and $t_1,\ldots,t_d$) such that, for every
$n\geq 1$:
\[
\sup\left|
E\left[\varphi
\left(\frac{Z_n(t_i)-Z_n(t_{i-1})}{\sqrt{t_i-t_{i-1}}}\right)_{1\leq
i \leq d}\right]-
E\big[\varphi(N)\big]\right|
\leq c\times
\left\{\begin{array}{lll}
\frac{1}{\sqrt{n}}&\,\,\mbox{if $H\in (0,\frac12]$}\\
\\
n^{2H-\frac32}&\,\,\mbox{if $H\in [\frac12,\frac34)$}\\
\\
\frac{1}{\sqrt{\log n}}&\,\,\mbox{if $H=\frac34$}
\end{array}\right.
\]
where the supremum is taken over all $C^2$-function $\varphi:\R^d\to\R$ such that $\|\varphi''\|_\infty\leq 1$.
\end{thm}
{\it Proof}.
We only make the proof for $H<3/4$, the proof for $H=3/4$ being similar.
Fix $d\geq 1$ and $t_0=0<t_1<\ldots<t_d$. In the sequel, $c$ will denote a constant independent of $n$,
which can differ from one line to another.
First, see e.g. the proof of Theorem \ref{BM}, observe that
$$\frac{Z_n(t_i)-Z_n(t_{i-1})}{\sqrt{t_i-t_{i-1}}}=I_2(f_i^{(n)})$$
with
$$
f_n^{(i)}=\frac{1}{\sigma_n\sqrt{t_i-t_{i-1}}
}
\sum_{k=\lfloor nt_{i-1}\rfloor}^{\lfloor nt_i\rfloor -1}
{\bf 1}_{[k,k+1]}^{\otimes 2}.
$$
In the
proof of Theorem \ref{BM}, it is shown that, for any fixed $i\in\{1,\ldots,d\}$ and $r\in\{1,\ldots,q_i-1\}$:
\begin{equation}\label{bound1}
\|f_n^{(i)}\otimes_1 f_n^{(i)}\|_{\HH^{\otimes 2}}\leq c\times
\left\{\begin{array}{lll}
\frac{1}{\sqrt{n}}&\,\,\mbox{if $H\in (0,\frac12]$}\\
\\
n^{2H-\frac32}&\,\,\mbox{if $H\in [\frac12,\frac34)$}\\
\end{array}\right..
\end{equation}
Moreover, when $1\leq i<j\leq d$, we have, with $\rho$ defined in (\ref{rho}),
\begin{eqnarray}\notag
&&\big|\langle f_n^{(i)},f_n^{(j)}\rangle_{\HH^{\otimes 2}}\big|\\
&=&\left|
\frac{1}{\sigma^2_n
\sqrt{t_i-t_{i-1}}\sqrt{t_j-t_{j-1}}
}
\sum_{k=\lfloor nt_{i-1}\rfloor}^{\lfloor nt_i\rfloor -1}
\sum_{l=\lfloor nt_{j-1}\rfloor}^{\lfloor nt_j\rfloor -1}
\rho^2(l-k)\right|\notag\\
&=& \frac{c}{\sigma^2_n}\left|
\sum_{|r|=\lfloor nt_{j-1}\rfloor-\lfloor nt_{i}\rfloor+1}
^{\lfloor nt_{j}\rfloor-\lfloor nt_{i-1}\rfloor-1}
\big[ ( \lfloor nt_{j}\rfloor -1 -r)\wedge( \lfloor nt_{i}\rfloor -1)
-( \lfloor nt_{j-1}\rfloor  -r)\vee( \lfloor nt_{i-1}\rfloor )
\big]
\rho^2(r)\right|\notag\\
&\leq& c\,\,\frac{\lfloor nt_{i}\rfloor -\lfloor nt_{i-1}\rfloor -1}{\sigma_n^2}
\sum_{|r|\geq\lfloor nt_{j-1}\rfloor-\lfloor nt_{i}\rfloor+1}
\rho^2(r)
=O\big(n^{4H-3}\big),\quad\mbox{as $n\to\infty$},\label{bound2}
\end{eqnarray}
the last equality coming from (\ref{sigman}) and
$$
\sum_{|r|\geq N}\rho^2(r)=O(\sum_{|r|\geq N}|r|^{4H-4})=O(N^{4H-3}),\quad\mbox{as $N\to\infty$}.
$$

Finally, by combining (\ref{bound1}), (\ref{bound2}), Corollary \ref{prop:Chaos} and
Proposition \ref{lm-control}, we obtain
the desired conclusion.
\fin

\end{document}